\documentclass[oneside,10pt]{article}          
\usepackage[b5paper]{geometry}	    
\usepackage{amsfonts,amsmath,latexsym,amssymb} 
\usepackage{theorem}                
\usepackage{mathrsfs,upref}         
\usepackage{mathptmx}		    
\usepackage{fdc}	            
\theoremstyle{definition}

\begin{document}

\title[Existence of solutions for a hybrid nonlinear generalized fractional pantograph equation]{About the existence of solutions for a hybrid nonlinear generalized fractional pantograph equation}


\author{Karimov E.T., L\'{o}pez B. and Sadarangani K.}

\address{Karimov E.T., Department of Mathematics and Statistics, Al-Khoud 123, Muscat, Oman,\\
\email{erkinjon@gmail.com}}

\address{L\'{o}pez B., Departamento de Matematicas, Universidad de Las Palmas de Gran Canaria, Campus de Tafira Baja, 35017 Las Palmas de Gran Canaria, Spain,\\
\email{blopez@dma.ulpgc.es}}

\address{Sadarangani K., Departamento de Matematicas, Universidad de Las Palmas de Gran Canaria, Campus de Tafira Baja, 35017 Las Palmas de Gran Canaria, Spain\\
\email{ksadaran@dma.ulpgc.es}}

\CorrespondingAuthor{Karimov E.T.}


\date{19.10.2015}                               




\begin{abstract}
        The main purpose of this paper is to study the existence of solutions for the following hybrid nonlinear fractional pantograph equation
$$
\left\{
\begin{aligned}
&D_{0+}^\alpha \left[\frac{x(t)}{f(t,x(t),x(\varphi(t)))}\right]=g(t,x(t),x(\rho(t))),\,\,0<t<1\\
&x(0)=0,
\end{aligned}
\right.
$$
where $\alpha\in (0,1)$, $\varphi$ and $\rho$ are functions from $[0,1]$ into itself and $D_{0+}^\alpha$ denotes the Riemann-Liouville fractional derivative. The main tool of our study is a generalization of Darbo's fixed point theorem associated to measures of non-compactness. Also, we present an example illustrating our results.
\end{abstract}

\maketitle



\section{Introduction}

        Fractional differential equations have recently been studied by a lot of number of researchers due to the fact that  they are valuable tools in the mathematical modelling of many phenomena appearing in the fields of physics, chemistry, aerodynamics, economics, control theory, signal and image processing, etc... For details, see, for example, [1-4] and the references therein.

In the literature, differential equations with proportional delays  are usually referred to as pantograph equations.  The name pantograph originated from the work [5] on the collection of current by the pantograph head of an electric locomotive. There are a great number of papers devoted to the qualitative properties  and numerical solutions of these equations (see, for example, [6-10]).

The pantograph equation has the form
$$
\left\{
\begin{aligned}
&y'(t)=ay(t)+by(\lambda t),\,\,0\leq t\leq T\\
&y(0)=y_0,
\end{aligned}
\right.
$$
where $0<\lambda<1$.

Recently, in [11], the authors considered the fractional version of the pantograph equation, namely
$$
\left\{
\begin{aligned}
&D_{0+}^\alpha u(t)=g(t,u(t),y(\lambda t)),\,\,0\leq t\leq T\\
&u(0)=u_0,
\end{aligned}
\right.
$$
where $\alpha, \lambda\in (0,1)$ and $D_{0+}^\alpha$ denotes the Riemann-Liouville fractional derivative. The main tool used in this study was the Banach contraction principle.

On the other hand, the following hybrid differential equation of first order
$$
\left\{
\begin{aligned}
&\frac{d}{dt} \left[\frac{x(t)}{f(t,x(t))}\right]=g(t,x(t)),\,\,0\leq t<T,\\
&x(t_0)=x_0
\end{aligned}
\right.
$$
was studied in [12] under the assumptions $f\in C\left([0,T)\times \mathbb{R}, \mathbb{R}\setminus\left\{0\right\}\right)$ and
$g\in C\left([0,T)\times \mathbb{R}, \mathbb{R}\right)$.

In [13], the authors discussed the fractional version of the last equation, i.e.,
$$
\left\{
\begin{aligned}
&D_{0+}^\alpha \left[\frac{x(t)}{f(t,x(t))}\right]=g(t,x(t)),\,\,0\leq t\leq T\\
&x(0)=0,
\end{aligned}
\right.
$$
where $\alpha\in (0,1)$, $f\in C\left([0,T]\times \mathbb{R}, \mathbb{R}\setminus\left\{0\right\}\right)$ and
$g\in C\left([0,T]\times \mathbb{R}, \mathbb{R}\right)$, being a fixed point theorem in Banach algebras the main tool used by the authors.

Recently, in [14] the authors studied the following hybrid fractional pantograph equation
$$
\left\{
\begin{aligned}
&D_{0+}^\alpha \left[\frac{x(t)}{f(t,x(t),x(\mu t))}\right]=g(t,x(t),x(\sigma t)),\,\,0<t<1\\
&x(0)=0,
\end{aligned}
\right.
$$
where $\alpha, \mu, \sigma \in (0,1)$,  $f\in C\left([0,1]\times \mathbb{R}\times \mathbb{R}, \mathbb{R}\setminus\left\{0\right\}\right)$ and
$g\in C\left([0,1]\times \mathbb{R}\times \mathbb{R}, \mathbb{R}\right)$, being a generalization of Darbo's fixed point theorem by using comparison functions the main tool used in the paper.

We also note work [15], where using a measure of non-compactness argument combined with the generalized version of Darbo's theorem, authors provide sufficient conditions for the existence of at least one solution to the functional equation
$$
x(t)=F\left(t,x(a(t)), \frac{f(t,x(b(t)))}{\Gamma_q(\alpha)}\int\limits_0^t (t-qs)^{(\alpha-1)u(s,x(s))d_qs}\right),\,\,t\in I,
$$
where $\alpha>1, \,q\in (0,1), \,I=[0,1],\,f,u:\,[0,1]\times \mathbb{R}\rightarrow \mathbb{R},\,a,b:\, I\rightarrow I$ and $F:\, I\times \mathbb{R}\times \mathbb{R}\rightarrow \mathbb{R}$.

In this paper, we study the following hybrid generalized fractional pantograph equation
\begin{equation}\label{eq1}
\left\{
\begin{aligned}
&D_{0+}^\alpha \left[\frac{x(t)}{f(t,x(t),x(\varphi(t)))}\right]=g(t,x(t),x(\rho(t))),\,\,0<t<1\\
&x(0)=0,
\end{aligned}
\right.
 \end{equation}
where $\alpha\in (0,1)$, $\varphi,\,\rho:\, [0,1]\rightarrow [0,1]$ are given functions.

The main tools in our study are a fixed point theorem for the product of two operators and a generalization of Darbo's fixed point theorem.
\section{Basic facts about fractional calculus}

In this section, we recall some definitions and basic results about fractional calculus. These results appear in [1].

\textbf{Definition 1.} The Riemann-Liouville fractional derivative of order $\alpha>0$ of a continuous function $f: (0,\infty)\rightarrow \mathbb{R}$ is given by
$$
D_{0+}^\alpha f(t)=\frac{1}{\Gamma(n-\alpha)} \left(\frac{d}{dt}\right)^{(n)}\int\limits_0^t \frac{f(s)}{(t-s)^{\alpha-n+1}}ds ,
$$
where $n=[\alpha]+1$, $[\alpha]$ denotes the integer part of $\alpha$ and $\Gamma$ denotes the classical gamma function, provided that the right side is point-wise defined on $(0,\infty))$.

\textbf{Definition 2.} The Riemann-Liouville fractional integral of order $\alpha>0$ of a continuous function $f: (0,\infty)\rightarrow \mathbb{R}$ is defined by
$$
I_{0+}^\alpha f(t)=\frac{1}{\Gamma(\alpha)} \int\limits_0^t {(t-s)^{\alpha-1}}{f(s)}ds ,
$$
provided that the right side is point-wise defined on $(0,\infty))$.

\textbf{Lemma 1.} Suppose that $f\in L^1(0,1)$ and $0<\alpha<1$. Then

(i) $D_{0+}^\alpha I_{0+}^\alpha f(t)=f(t)$.

(ii) $I_{0+}^\alpha D_{0+}^\alpha f(t)=f(t)- \frac{I_{0+}^\alpha f(t)|_{t=0}}{\Gamma(\alpha)}t^{\alpha-1}$ a.e. in $(0,1)$.

\bigskip

Following the same argument which appears in the proof of Lema 2.4 of [14], the following lemma can be proved.

\textbf{Lemma 2.} Let $0<\alpha<1$ and suppose that $f\in C\left([0,1]\times \mathbb{R}\times \mathbb{R}, \mathbb{R}\setminus \{0\}\right)$ and $h\in C[0,1]$. Then the unique solution of the fractional hybrid initial value problem
$$
\left\{
\begin{aligned}
&D_{0+}^\alpha \left[\frac{x(t)}{f(t,x(t),x(\varphi(t)))}\right]=h(t),\,\,0<t<1\\
&x(0)=0,
\end{aligned}
\right.
$$
where $\varphi:\, [0,1]\rightarrow [0,1]$ is a continuous function, is given by
$$
x(t)=\frac{f\left(t,x(t),x(\varphi(t))\right)}{\Gamma(\alpha)}\int\limits_0^t \frac{h(s)}{(t-s)^{1-\alpha}}ds,\,\,t\in [0,1].
$$

\section{Background about measures of non-compactness}

In this section we present some facts and basic results about measures of non-compactness which will be used later.

Assume that $E$ is a real Banach space with nor $||\cdot||$ and the zero element $0$. By $B(x,r)$ we denote the closed ball in $E$ centered at $x$ with the radius $r$. By $B_r$ we denote the ball $B(0,r)$. If $X$ is non-empty subset of $E$, then $\overline{X}$ and $Conv X$ denote the closure and the closed convex closure of $X$, respectively. When $X$ is a bounded subset, $diam X$ denotes the diameter of $X$ and $||X||$ the quantity given by $||X||=sup\{||x||:\, x\in X\}$. Further, by $\mathfrak{M}_E$ we denote the family of the non-empty and bounded subsets of $E$ and by $\mathfrak{N}_E$ its subfamily consisting of the relatively compact subsets.

In the paper, we accept the following definition of measure of non-compactness which appears in [16].

\textbf{Definition 3.} A mapping $\mu:\, \mathfrak{M}_E\rightarrow \mathbb{R}=[0,\infty)$ will be called a measure of non-compactness in $E$ if it satisfies the following conditions:

$1^\circ.\,\,$ The family $Ker \mu =\{X\in \mathfrak{M}_E;\,\mu(X)=0\}$ is non-empty and $Ker \mu\in \mathfrak{N}_E$.

$2^\circ.\,\,$ $X\subset Y \Rightarrow \mu(X)\leq \mu(Y)$.

$3^\circ.\,\,$ $\mu(\overline{X})=\mu (Conv X)=\mu(X)$.

$4^\circ.\,\,$ $\mu (\lambda X+(1-\lambda)Y)\leq \lambda\mu(X)+(1-\lambda)\mu(Y)$ for $\lambda\in[0,1]$.

$5^\circ.\,\,$ If $(X_n)$ is a sequence of closed subsets of $\mathfrak{M}_E$ such that $X_{n+1}\subset X_n$ for $n=1,2,3,...,$ and $\lim\limits_{n\rightarrow \infty}\mu(X_n)=0$ then $X_\infty=\bigcap\limits_{n=1}^\infty X_n\neq \emptyset$.

The family $Ker \mu$ appearing in $1^\circ$  is called the kernel of the measure of non-compactness $\mu$. Notice that the set $X_\infty$ appearing in $5^\circ$ belongs to $Ker \mu$. Indeed, since $\mu(X_\infty)\subset \mu(X_n)$ for any $n=1,2,3,...,$ it follows that $\mu(X_\infty)\leq\lim\limits_{n\rightarrow \infty} \mu(X_n)=0$.

In connection with measures of non-compactness, in [17] Darbo proved the following fixed point theorem.

\textbf{Theorem 1.} Let $\Omega$ be a nonempty, bounded, closed and convex subset of a Banach space $E$ and let $T:\,\Omega\rightarrow \Omega$ be a continuous mapping. Suppose that there exists $k\in [0,1)$ such that
$$
\mu(TX)\leq k \mu(X),
$$
for any non-empty subset $X$ of $\Omega$, where $\mu$ is a measure of non-compactness in $E$. Then $T$ has a fixed point in $\Omega$.

Recently, some generalizations of Theorem 1 have appeared in the literature (see [18-21], for example). The following generalization of Darbo's fixed point theorem appears in [21] and it is the version in the context of measures of non-compactness of a recent result about fixed point theorem which appears in [22]. For the paper is self-contained, we present this result. Previously, we need to introduce the class $\mathfrak{F}$ of functions. By $\mathfrak{F}$ we denote the class of functions $\varphi:\, (0,\infty)\rightarrow (1,\infty)$ satisfying the following condition:

For any sequence $(t_n)\subset (0,\infty)$
$$
\lim\limits_{n\rightarrow\infty}\varphi(t_n)=1 \Leftrightarrow \lim\limits_{n\rightarrow\infty}t_n=0.
$$
Examples of functions belonging to the class of $\mathfrak{F}$  are $\varphi (t)=e^{\sqrt{t}}$, $\varphi(t)=2-\frac{2}{\pi}\arctan(\frac{1}{t^\alpha})$ with $0<\alpha<1$ [21] and $\varphi(t)=(1+t^2)^\beta$ with $\beta>0$.

\textbf{Theorem 2.} Let $\Omega$ be a nonempty, bounded, closed and convex subset of a Banach space $E$ and let $T:\,\Omega\rightarrow \Omega$ be a continuous mapping. Suppose that there exist $\varphi\in \mathbb{F}$ and $k\in (0,1)$ such that, for any nonempty subset $X$ of $\Omega$ with $\mu(TX)>0$,
$$
\varphi\left(\mu(TX)\right)\leq \left(\varphi\left(\mu(X)\right)\right)^k,
$$
 where $\mu$ is a measure of non-compactness in $E$. Then $T$ has a fixed point in $\Omega$.

Next, we introduce the following concept which appears in [23] which will be important for our purposes.

\textbf{Definition 4.} Let $E$ be a Banach algebra. A measure of non-compactness $\mu$ in $E$ said to satisfy condition (m) if it satisfies the following condition:
$$
\mu(XY)\leq ||X||\mu(Y)+||Y||\mu(X)
$$
for any $X, Y\in \mathfrak{M}_E$, where $XY=\{xy:\, x\in X, y\in Y\}$.

In the paper we work in the space $C[0,1]$ of the real functions defined and continuous on $[0,1]$, with the usual supremum norm given by
$||x||=\sup \{|x(t)|:\, t\in [0,1]\}$ for $x\in C[0,1]$. Notice that $\left(C[0,1], ||\cdot||\right)$ is a Banach algebra, where the multiplication is defined as the usual product of real functions.

Next, we present measure of non-compactness in $C[0,1]$ which will be used in our study. Fix a set $X\in \mathfrak{M}_{C[0,1]}$ and $\varepsilon>0$. For $x\in X$, by $w(x,\varepsilon)$ we denote the modulus of continuity of $x$, i.e.,
$$
w(x,\varepsilon)=\sup \{|x(t)-x(s)|:\, t,s\in[0,1], |t-s|\leq \varepsilon\}.
$$
Further, put
$$
w(X,\varepsilon)=\sup \{w(x,\varepsilon):\, x\in X\}
$$
and
$$
w_0(X)=\lim\limits_{\varepsilon\rightarrow 0}w(X,\varepsilon).
$$
In [16], it is proved that $w_0$ is a measure of non-compactness in $C[0,1]$.

\section{Main result}
Problem (\ref{eq1}) will be studied under the following assumptions:

\textbf{(H1):} $f\in C\left([0,1]\times \mathbb{R}\times \mathbb{R}, \mathbb{R}\setminus \{0\}\right)$ and $g\in C\left([0,1]\times \mathbb{R}\times \mathbb{R}, \mathbb{R}\right)$.

\textbf{(H2):} $\varphi,\rho: [0,1]\rightarrow [0,1]$ are continuous function.

\textbf{(H3):} The functions $f$ and $g$ satisfy
$$
\left|f(t,x_1,y_1)-f(t,x_2,y_2)\right|\leq \left(max (|x_1-x_2|, |y_1-y_2|)+1\right)^k-1
$$
and
$$
\left|g(t,x_1,y_1)-g(t,x_2,y_2)\right|\leq \left(max (|x_1-x_2|, |y_1-y_2|)+1\right)^r-1,
$$
respectively, for any $t\in [0,1]$ and $x_1, x_2, y_1, y_2\in \mathbb{R}$, where $k,r\in (0,1)$.

Notice that assumption (H1) gives us the existence of two nonnegative constants $K_1$ and $K_2$ such that
$$
K_1=sup\left\{|f(t,0,0)|: t\in [0,1]\right\},\,\,\,K_2=sup\left\{|g(t,0,0)|: t\in [0,1]\right\}.
$$

\textbf{(H4):} There exists $r_0>0$ such that
$$
\left((r_0+1)^k-1+K_1\right)\left((r_0+1)^r-1+K_2\right)\leq r_0\Gamma(\alpha+1)
$$
and
$$
(r_0+1)^r-1+K_2\leq \Gamma(\alpha+1).
$$

\textbf{Theorem 3.}
Under assumptions (H1)-(H4), Problem (\ref{eq1}) has at least one solution in $C[0,1]$.

\textbf{Proof:}

Consider the operator $\mathbb{T}$ defined on $C[0,1]$ by
$$
(\mathbb{T}x)(t)=\frac{f\left(t,x(t),x(\varphi(t))\right)}{\Gamma(\alpha)}\int\limits_0^t\frac{g(s,x(s),x(\rho(s)))}{(t-s)^{1-\alpha}}ds,
$$
for $x\in C[0,1]$ and $t\in [0,1]$

In virtue of Lemma 2.4, a fixed point of $\mathbb{T}$ gives us the desired result.

Let $\mathbb{F}$ and $\mathbb{G}$ be the operators defined on $C[0,1]$ by
$$
(\mathbb{F}x)(t)=f(t,x(t),x(\varphi(t)))
$$
and
$$
(\mathbb{G}x)(t)=\frac{1}{\Gamma(\alpha)}\int\limits_0^t\frac{g(s,x(s),x(\rho(s)))}{(t-s)^{1-\alpha}}ds,
$$
for any $x\in C[0,1]$ and $t\in [0,1]$.

Then $\mathbb{T}x=(\mathbb{F}x)\cdot(\mathbb{G}x)$ for any $x\in C[0,1]$.

For a better readability, we divide the proof in several steps.

\textbf{Step 1:} $\mathbb{T}$ applies $C[0,1]$ into itself.

In fact, since the product of continuous functions is a continuous function, it is sufficient to prove that $\mathbb{F}x, \mathbb{G}x\in C[0,1]$ for any $x\in C[0,1]$.

It is clear , by (H1) and (H2), that if $x\in C[0,1]$ then $\mathbb{F}x\in C[0,1]$. Next, we will prove that if
$x\in C[0,1]$ then $\mathbb{G}x\in C[0,1]$.

To do this, we fix $t_0\in [0,1]$ and let $(t_n)$ be a sequence in $[0,1]$ such that $t_n\rightarrow t_0$.

In fact, without loss of generality, we can suppose that $t_n>t_0$. Then,
$$
\begin{aligned}
&\left|(\mathbb{G}x)(t_n)-(\mathbb{G}x)(t_0)\right|=\frac{1}{\Gamma(\alpha)}\left|\int\limits_0^{t_n}\frac{g(s,x(s),x(\rho(s)))}{(t_n-s)^{1-\alpha}}ds-\int\limits_0^{t_0}\frac{g(s,x(s),x(\rho(s)))}{(t_0-s)^{1-\alpha}}ds\right|\leq\\
&\leq \frac{1}{\Gamma(\alpha)}\left|\int\limits_0^{t_n}\frac{g(s,x(s),x(\rho(s)))}{(t_n-s)^{1-\alpha}}ds-\int\limits_0^{t_n}\frac{g(s,x(s),x(\rho(s)))}{(t_0-s)^{1-\alpha}}ds\right|\\
&+\frac{1}{\Gamma(\alpha)}\left|\int\limits_0^{t_n}\frac{g(s,x(s),x(\rho(s)))}{(t_0-s)^{1-\alpha}}ds-\int\limits_0^{t_0}\frac{g(s,x(s),x(\rho (s)))}{(t_0-s)^{1-\alpha}}ds\right|\\
&\leq \frac{1}{\Gamma(\alpha)}\int\limits_0^{t_n}\left|(t_n-s)^{\alpha-1}-(t_0-s)^{\alpha-1}\right|\left|g(s,x(s),x(\rho(s)))\right|ds\\
&+\frac{1}{\Gamma(\alpha)}\int\limits_{t_0}^{t_n}\left|t_0-s\right|^{\alpha-1}\left|g(s,x(s),x(\rho(s)))\right|ds.\\
\end{aligned}
$$
Since $g\in C\left([0,1]\times \mathbb{R}\times \mathbb{R}, \mathbb{R}\right)$, $g$ will be bounded on the compact $[0,1]\times[-||x||, ||x||]\times [-||x||, \left\|x\right\|]$ and we put
$$
L=sup\left\{|g(s,x,y)|:\, s\in[0,1],\, x,y\in [-||x||, \left\|x\right\|]\right\}.
$$
From the last estimate, we get
$$
\left|(\mathbb{G}x)(t_n)-(\mathbb{G}x)(t_0)\right|\leq \frac{L}{\Gamma(\alpha)}\int\limits_0^{t_n}\left|(t_n-s)^{\alpha-1}-(t_0-s)^{\alpha-1}\right|ds+\frac{L}{\Gamma(\alpha)}\int\limits_{t_0}^{t_n}|t_0-s|^{\alpha-1}ds.
$$
Taking into account that $0<\alpha<1$ and $t_n>t_0$, we infer
$$
\begin{aligned}
&\left|(\mathbb{G}x)(t_n)-(\mathbb{G}x)(t_0)\right|\leq \\ &\frac{L}{\Gamma(\alpha)}\left[\int\limits_0^{t_0}\left|(t_n-s)^{\alpha-1}-(t_0-s)^{\alpha-1}\right|ds+\int\limits_{t_0}^{t_n}\left|(t_n-s)^{\alpha-1}-(t_0-s)^{\alpha-1}\right|ds\right]+
\frac{L}{\Gamma(\alpha)}\int\limits_{t_0}^{t_n}|s-t_0|^{\alpha-1}ds\\
&=\frac{L}{\Gamma(\alpha)}\left[\int\limits_0^{t_0}\left[(t_0-s)^{\alpha-1}-(t_n-s)^{\alpha-1}\right]ds+
\int\limits_{t_0}^{t_n}\frac{ds}{(t_n-s)^{1-\alpha}}+\int\limits_{t_0}^{t_n}\frac{ds}{(s-t_0)^{1-\alpha}}\right]\\
&+\frac{L}{\Gamma(\alpha)}\int\limits_{t_0}^{t_n}\frac{ds}{(s-t_0)^{1-\alpha}}
<\frac{4L}{\Gamma(\alpha+1)}(t_n-t_0)^\alpha,
\end{aligned}
$$
where we have used the fact that $t_0^\alpha-t_n^\alpha<0$.

From the last estimate, we deduce that $(\mathbb{G}x)(t_n)\rightarrow (\mathbb{G}x)(t_0)$ when $n\rightarrow \infty$. This proves that if $x\in C[0,1]$.

\textbf{Step 2:} An estimate of $||\mathbb{T}x||$ for $x\in C[0,1]$.

Fix $x\in C[0,1]$ and $t\in C[0,1]$. In view of assumptions, we have
$$
\begin{aligned}
&|(\mathbb{T}x)(t)|=|(\mathbb{F}x)(t)|\cdot (\mathbb{G}x)(t)=|f(t,x(t),x(\varphi(t)))|\cdot\left|\frac{1}{\Gamma(\alpha)}\int\limits_0^t\frac{g(s,x(s),x(\rho(s)))}{(t-s)^{1-\alpha}}ds\right|\\
&\leq \left[|f(t,x(t),x(\varphi(t)))-f(t,0,0)|+|f(t,0,0)|\right]\times\\ &\frac{1}{\Gamma(\alpha)}\left|\int\limits_0^t\frac{g(s,x(s),x(\rho(s)))-g(s,0,0)}{(t-s)^{1-\alpha}}ds+\int\limits_0^t\frac{g(s,0,0)}{(t-s)^{1-\alpha}}ds\right|\\
&\leq \frac{1}{\Gamma(\alpha)}\left[(max(|x(t)|,|x(\rho(t))|)+1)^k-1+K_1\right]\times\\
&\left[\int\limits_0^t\frac{(max(|x(s)|, |x(\varphi(s))|)+1)^r-1}{(t-s)^{1-\alpha}}ds+K_2\int\limits_0^t\frac{ds}{(t-s)^{1-\alpha}}\right]\\
&\leq \frac{1}{\Gamma(\alpha)}\left[(max(||x||, ||x||)+1)^k-1+K_1\right]\left[(max(||x||, ||x||)+1)^r+K_2\right]\int\limits_0^t\frac{ds}{(t-s)^{1-\alpha}}\\
&\leq \frac{1}{\Gamma(\alpha)}\left[(||x||+1)^k-1+K_1\right]\left[(||x||+1)^r-1+K_2\right]\frac{t^\alpha}{\alpha}
\\
&\leq \frac{1}{\Gamma(\alpha+1)}\left[(||x||+1)^k-1+K_1\right]\left[(||x||+1)^r-1+K_2\right].
\end{aligned}
$$
Therefore,
$$
||\mathbb{T}x||\leq \frac{1}{\Gamma(\alpha+1)}\left[(||x||+1)^k-1+K_1\right]\left[(||x||+1)^r-1+K_2\right].
$$
By assumption (H4), we infer that the operator $\mathbb{T}$ applies $B_{r_0}$ into itself. Moreover, from the last estimates, it follows that
$$
||\mathbb{F} B_{r_0}||\leq (r_0+1)^k-1+K_1
$$
and
$$
||\mathbb{G} B_{r_0}||\leq \frac{1}{\Gamma(\alpha+1)}\left[(r_0+1)^r-1+K_2\right].
$$

\textbf{Step 3:} The operators $\mathbb{F}$ and $\mathbb{G}$ are continuous on the ball $B_{r_0}$.

In fact, firstly we prove that $\mathbb{F}$ is continuous on $B_{r_0}$. To do this, we fix $\varepsilon>0$ and we take $x,y\in B_{r_0}$ with $||x-y||\leq \varepsilon$. Then, for $t\in [0,1]$, we have
$$
\begin{aligned}
&|(\mathbb{F}x)(t)-(\mathbb{F}y)(t)|=|f(t,x(t),x(\varphi(t)))-f(t,y(t),y(\varphi(t)))|\\
&\leq (max(|x(t)-y(t)|, |x(\varphi(t))-y(\varphi(t))|)+1)^k-1\\
&\leq (max(||x-y||, ||x-y||+1)+1)^k-1=(||x-y||+1)^k-1\leq (\varepsilon+1)^k-1\\
\end{aligned}
$$
and, since $(\varepsilon+1)^k-1\rightarrow 0$ when $\varepsilon\rightarrow 0$, we have proved that $\mathbb{F}$ is continuous in $B_{r_0}$.

Next, we prove that $\mathbb{G}$ is continuous in $B_{r_0}$. In order to do this, we  fix $\varepsilon>0$ and we take $x,y\in B_{r_0}$ with $||x-y||\leq \varepsilon$. Then, for $t\in [0,1]$, we get
$$
\begin{aligned}
&|(\mathbb{G}x)(t)-(\mathbb{G}y)(t)|=\frac{1}{\Gamma (\alpha)}\left|\int\limits_0^t \frac{g(s,x(s),x(\rho(s)))}{(t-s)^{1-\alpha}}ds-
\int\limits_0^t \frac{g(s,y(s),y(\rho(s)))}{(t-s)^{1-\alpha}}ds\right|\\
&\leq \frac{1}{\Gamma(\alpha)}\int\limits_0^t \frac{|g(s,x(s),x(\rho(s)))-g(s,y(s),y(\rho(s)))|}{(t-s)^{1-\alpha}}ds\\
&\leq \frac{1}{\Gamma(\alpha)}\int\limits_0^t \frac{(max(|x(s)-y(s)|, |x(\rho(s))-y(\rho(s))|)+1)^r-1}{(t-s)^{1-\alpha}}ds\\
&=\frac{1}{\Gamma(\alpha)}\int\limits_0^t \frac{(max(||x-y||, ||x-y||)+1)^r-1}{(t-s)^{1-\alpha}}ds
=\frac{\left[(||x-y||+1)^r-1\right]}{\Gamma(\alpha)}\int\limits_0^t(t-s)^{\alpha-1}ds\\
&=\frac{\left[(||x-y||+1)^r-1\right]}{\Gamma(\alpha)}
\cdot\frac{t^\alpha}{\alpha}\leq \frac{(\varepsilon+1)^r-1}{\Gamma(\alpha+1)}
\end{aligned}
$$
and, as $\frac{(\varepsilon+1)^r-1}{\Gamma(\alpha+1)}\rightarrow 0$ when $\varepsilon \rightarrow 0$, we have proved that $\mathbb{G}$ is continuous on $B_{r_0}$. Consequently, since $\mathbb{T}=\mathbb{F}\cdot \mathbb{G}$, it follows that $\mathbb{T}$ is continuous on $B_{r_0}$.

\textbf{Step 4:} Estimates of $\omega_0(\mathbb{F}X)$ and $\omega_0(\mathbb{G}X)$ for $\emptyset\neq X\subset B_{r_0}$.

Firstly, we estimate $\omega_0(\mathbb{F}X)$. For $\varepsilon>0$ given, since $\varphi:\, [0,1]\rightarrow [0,1]$ is uniformly continuous, we can find $\delta>0$ (which can be taken with $\delta<\varepsilon$) such that, for $|t_1-t_2|<\delta$ we have $|\varphi(t_1)-\varphi(t_2)|<\varepsilon$.

Now, we take $x\in X$ and $t_1, t_2\in [0,1]$ with $|t_1-t_2|\leq \delta<\varepsilon$. Then
$$
\begin{aligned}
&|(\mathbb{F}x)(t_1)-(\mathbb{F}x)(t_2)|=|f(t_1,x(t_1),x(\varphi(t_1)))-f(t_2,x(t_2),x(\varphi(t_2)))|\\
&\leq |f(t_1,x(t_1),x(\varphi(t_1)))-f(t_1,x(t_2),x(\varphi(t_2)))|+|f(t_1,x(t_2),x(\varphi(t_2)))-f(t_2,x(t_2),x(\varphi(t_2)))|\\
&\leq \left[(max(|x(t_1)-x(t_2)|, |x(\varphi(t_1))-x(\varphi(t_2))|)+1)^k-1\right]+w(f,\varepsilon)\\
&\leq \left[(w(X,\varepsilon)+1)^k-1\right]+w(f,\varepsilon),
\end{aligned}
$$
where $w(f,\varepsilon)$ denotes the quantity
$$
w(f,\varepsilon)=sup\left\{|f(t_1,x,y)-f(t_2,x,y)|:\, t_1,t_2\in [0,1],\, |t_1-t_2|\leq \delta,\, x,y\in [-r_0, r_0]\right\}.
$$
Therefore,
$$
w(\mathbb{F}X,\delta)\leq \left[(w(X,\varepsilon)+1)^k-1\right]+w(f,\varepsilon).
$$
Since $f(t,x,y)$ is uniformly continuous on the compact $[0,1]\times [-r_0, r_0]\times [-r_0, r_0]$, $w(f,\varepsilon)\rightarrow 0$ when $\varepsilon\rightarrow 0$, and, consequently, from the last inequality, we infer
$$
\omega_0(\mathbb{F}X)\leq \left[(\omega_0(X)+1)^k-1\right]
$$
Next, we estimate $\omega_0(\mathbb{G}X)$. Fix $\varepsilon>0$, and we take $x\in X$ and $t_1, t_2\in [0,1]$ with $|t_1-t_2|\leq \varepsilon$. Without loss of generality, we can suppose that $t_1<t_2$. Then, we have
$$
\begin{aligned}
&|(\mathbb{G}x)(t_2)-(\mathbb{G}x)(t_1)|=\frac{1}{\Gamma (\alpha)}\left|\int\limits_0^{t_2} \frac{g(s,x(s),x(\rho(s)))}{(t_2-s)^{1-\alpha}}ds-
-\int\limits_0^{t_1} \frac{g(s,x(s),x(\rho(s)))}{(t_1-s)^{1-\alpha}}ds\right|\\
&\leq \frac{1}{\Gamma (\alpha)}\left[\int\limits_0^{t_1}|(t_2-s)^{\alpha-1}-(t_1-s)^{\alpha-1}||g(s,x(s),x(\rho(s)))|ds\right.\\
&\left.+\int\limits_{t_1}^{t_2}(t_2-s)^{\alpha-1}|g(s,x(s),x(\rho(s)))|ds\right]\\
&=\frac{1}{\Gamma (\alpha)}\left[\int\limits_0^{t_1}\left[(t_1-s)^{\alpha-1}-(t_2-s)^{\alpha-1}\right]|g(s,x(s),x(\rho(s)))|ds\right.\\
&\left.+\int\limits_{t_1}^{t_2}(t_2-s)^{\alpha-1}|g(s,x(s),x(\rho(s)))|ds\right].
\end{aligned}
$$
Since $g(t,x,y)$ is continuous on $[0,1]\times \mathbb{R}\times \mathbb{R}$, it is bounded on the compact subset $[0,1]\times [-r_0, r_0]\times [-r_0, r_0]$. Put $M=sup\left\{|g(t,x,y):\, t\in [0,1],\, x,y\in [-r_0, r_0]\right\}$. then, from the last estimate, we infer that
$$
\begin{aligned}
&\left|(\mathbb{G}x)(t_2)-(\mathbb{G}x)(t_1)\right|=\frac{M}{\Gamma (\alpha)}\left[\int\limits_0^{t_1}\left[(t_1-s)^{\alpha-1}-(t_2-s)^{\alpha-1}\right]ds+
\int\limits_{t_1}^{t_2}(t_2-s)^{\alpha-1}ds\right]\\
&\leq \frac{M}{\Gamma (\alpha+1)}\left[(t_2-t_1)^\alpha+t_1^\alpha-t_2^\alpha+(t_2-t_1)^\alpha\right]\leq \frac{2M}{\Gamma (\alpha+1)}(t_2-t_1)^\alpha\leq \frac{2M}{\Gamma(\alpha+1)}\varepsilon^\alpha,
\end{aligned}
$$
where we have used the fact that $t_1^\alpha-t_2^\alpha\leq 0$. Therefore,
$$
w(\mathbb{G}x,\varepsilon)\leq \frac{2M}{\Gamma(\alpha+1)}\varepsilon^\alpha
$$
and this gives us $\omega_0(\mathbb{G}X)=0$.

\textbf{Step 5:} An estimate of $\omega_0(\mathbb{T}X)$ for $\emptyset\neq X\subset B_{r_0}$.

Taking into account that
$$
\omega_0(XY)\leq ||X||\omega_0(Y)+||Y||\omega_0(X)
$$
from the estimates obtained in steps 2 and 4, we deduce
$$
\begin{aligned}
&\omega_0(\mathbb{T}X)=\omega_0(\mathbb{F}X\cdot\mathbb{G}X)\leq ||\mathbb{F}X||\omega_0(\mathbb{G}X)+||\mathbb{G}X||\omega_0(\mathbb{F}X)\\
&\leq ||\mathbb{F}B_{r_0}||\omega_0(\mathbb{G}X)+||\mathbb{G}B_{r_0}||\omega_0(\mathbb{F}X)\\
&\leq \frac{1}{\Gamma(\alpha+1)}\left[(r_0+1)^r-1+K_2\right]\left[(\omega_0(X)+1)^k-1\right].
\end{aligned}
$$
By assumption (H4), $\frac{1}{\Gamma(\alpha+1)}[(r_0+1)^r-1+K_2]\leq 1$ and, from the last estimate, we infer that
$$
\omega_0(\mathbb{T}X)\leq (\omega_0(X)+1)^k-1
$$
or, equivalently,
$$
\omega_0(\mathbb{T}X)+1\leq (\omega_0(X)+1)^k.
$$
Therefore, the contractive condition appearing in Theorem 1 is satisfied with $\varphi(t)=t+1$, where $\varphi\in \mathfrak{F}$. By Theorem 2, the operator $\mathbb{T}$ has at least one fixed point in $B_{r_0}$.

This completes the proof.

\section{Examples and comparison with other results}

In this section we present an example illustrating our results. Previously, we will need the following lemma:

\textbf{Lemma 3.} Let $\varphi:\, [0,\infty)\rightarrow [0,\infty)$ be the function defined by $\varphi(t)=(t+1)^k-1$ for $t\in [0,\infty)$, where $k\in (0,1)$. then, we have

(a) $\varphi$ is nondecreasing;

(b) $|\varphi(t)-\varphi(t')|\leq \varphi(|t-t'|)$ for any $t,t'\in [0,\infty)$.

\textbf{Proof:}

(a) It is clear since $\varphi'(t)=k(t+1)^{k-1}>0$ for $t\in [0,\infty)$.

(b) Since $\varphi''(t)=k(k-1)(t+1)^{k-2}<0$ for $t\in [0,\infty)$, $\varphi$ is a concave function. Moreover, $\varphi(0)=0$.

It is a well known fact that the concavity of $\varphi$ and $\varphi(0)=0$ imply the subadditivity of the function $\varphi$, i.e.
$\varphi(t+t')\leq \varphi(t)+\varphi(t')$ for any $t,t'\in [0,\infty)$.

In order to prove (b), we can suppose without loss of generality that $t<t'$. Then
$$
\varphi(t')=\varphi(t'-t+t)\leq \varphi(t'-t)+\varphi(t)
$$
and, consequently, $\varphi(t')-\varphi(t)\leq \varphi(t'-t)$ and this proves (b).

\textbf{Example 1.} Consider the following fractional hybrid problem
\begin{equation}\label{eq2}
\left\{
\begin{aligned}
&D_{0+}^{1/2} \left[\frac{x(t)}{\frac{1}{2}\left(\sqrt[4]{1+|x(t)|}+\sqrt[4]{1+\left|x\left(\frac{t}{1+t}\right)\right|}\right)}\right]\\
&=\frac{1}{\beta}
\left[\sqrt[3]{1+|x(t)|}+\sqrt[3]{1+\left|x\left(arctg t\right)\right|}\right],\,0<t<1\\
&x(0)=0,
\end{aligned}
\right.
\end{equation}
where $\alpha,\beta>0$.

Notice that Problem (2) is a particular case of Problem (1), with $\alpha=1/2, \varphi(t)=\frac{t}{1+t}$, $\rho(t)=arctg t,$ $f(t,x,y)=\frac{1}{\alpha}\left(\sqrt[4]{1+|x|}+\sqrt[4]{1+|y|}\right)$ and $g(t,x,y)=\frac{1}{\beta}\left(\sqrt[3]{1+|x|}+\sqrt[3]{1+|y|}\right)$.

It is clear that assumptions (H1) and (H2) of the Theorem 3 are satisfied. On the other hand, for any $t\in [0,1]$ and $x,y,x_1,y_1\in \mathbb{R}$, we have
$$
\begin{aligned}
&|f(t,x,y)-f(t,x_1,y_1)|=\left|\frac{1}{\alpha}\left(\sqrt[4]{1+|x|}+\sqrt[4]{1+|y|}\right)-\frac{1}{\alpha}\left(\sqrt[4]{1+|x_1|}+\sqrt[4]{1+|y_1|}\right)\right|\\
&\leq \frac{1}{\alpha}\left|\sqrt[4]{1+|x|}-\sqrt[4]{1+|x_1|}\right|+\frac{1}{\alpha}\left|\sqrt[4]{1+|y|}-\sqrt[4]{1+|y_1|}\right|\\
&\leq \frac{1}{\alpha}\left(\sqrt[4]{1+|x|}-1\right)-\left(\sqrt[4]{1+|x_1|}-1\right)+\frac{1}{\alpha}\left(\sqrt[4]{1+|y|}-1\right)-\left(\sqrt[4]{1+|y_1|}-1\right)\\
&\leq \frac{1}{\alpha}\left(\sqrt[4]{||x|-|x_1||+1}-1\right)+\frac{1}{\alpha}\left(\sqrt[4]{||y|-|y_1||+1}-1\right)\\
&\leq \frac{2}{\alpha}
\left(\sqrt[4]{max(|x-x_1|, |y-y_1|)+1}-1\right),
\end{aligned}
$$
where we have used Lemma 3.

By using a similar argument, we can prove
$$
|g(t,x,y)-g(t,x_1,y_1)|\leq \frac{2}{\beta}
\left(\sqrt[3]{max(|x-x_1|, |y-y_1|)+1}-1\right)
$$
for any $t\in [0,1]$ and $x,y,x_1,y_1\in \mathbb{R}$.

Therefore, assumption (H3) of the Theorem 3 is satisfied when $\alpha,\beta\geq 2$ with $k=1/4$ and $r=1/3$.

In our case, $K_1=sup\{|f(t,0,0)|:\, t\in[0,1]\}=\frac{2}{\alpha}$ and $K_2=sup\{|g(t,0,0)|:\, t\in[0,1]\}=\frac{2}{\beta}$.

The inequality appearing in (H4) of Theorem 3 has the expression
$$
\left((r_0+1)^{1/4}-1+\frac{2}{\alpha}\right)\left((r_0+1)^{1/3}-1+\frac{2}{\beta}\right)\leq r_0\Gamma(3/2).
$$
For $\alpha=4$ and $\beta=3$, we have
$$
\left((r_0+1)^{1/4}-\frac{1}{2}\right)\left((r_0+1)^{1/3}-\frac{12}{3}\right)\leq r_0\Gamma(3/2).
$$
This inequality is satisfied for $r_0=0.8$. Moreover,
$$
\left((r_0+1)^{1/3}-\frac{1}{3}\right)=\sqrt[3]{1.8}-1/3\cong 0.8832\leq \Gamma(3/2)=0.8862.
$$
Therefore, Problem (2) has at least one solution $x(t)\in C[0,1]$ for $\alpha=4$ and $\beta=3$ with $||x||\leq 0.8$.

\bigskip

An interesting question is the non-oscillatory character of the solutions of Problem (\ref{eq1}), i.e., that the solutions of Problem (\ref{eq1}) have a constant sign. Notice that if $f(t,x,y)$ and $g(t,x,y)$ have the same sign (this means that $f(t,x,y)>0$ and $g(t,x,y)\geq 0$ or $f(t,x,y)$ and $g(t,x,y)\leq 0$ for any $t\in [0,1]$ and $x,y\in \mathbb{R}$) and assumptions of Theorem 3 are satisfied, then the solutions $x(t)$ of Problem (1) are nonnegative due to the fact that these solutions satisfy the integral equation
$$
x(t)=\frac{f(t,x(t),x(\varphi(t)))}{\Gamma(\alpha)}\int\limits_0^t \frac{g(s,x(s),x(\rho(s)))}{(t-s)^{1-\alpha}}ds,\, 0\leq t\leq 1.
$$
In connection with the above-mentioned question, we have the following result.

\textbf{Proposition 1.} Under assumptions of Theorem 3 and suppose that $g(t,x,y)$ has constant sign and $g(t,x,y)\neq 0$ for $t\in[0,1]$ and $x,y\in \mathbb{R}$, we have that the solution $x(t)$ of Problem (1) obtained in Theorem 3 satisfies that $x(t)\neq 0$ for $0<t<1$.

\textbf{Proof:} Suppose the contrary case. Then we can find $t^*\in (0,1)$ with $x(t^*)=0$. Since $x(t)$ satisfies the last integral equation, we have
$$
0=x(t^*)=\frac{f(t^*,x(t^*),x(\varphi(t^*)))}{\Gamma(\alpha)}\int\limits_0^{t^*} \frac{g(s,x(s),x(\rho(s)))}{(t^*-s)^{1-\alpha}}ds.
$$
Taking into account that $f(t,x,y)\neq 0$ for any $t\in [0,1]$ and $x,y\in\mathbb{R}$ we infer that
$$
\int\limits_0^{t^*} (t^*-s)^{1-\alpha}g(s,x(s),x(\rho(s)))ds=0.
$$
Since $g(t,x,y)$ has constant sign and $(t^*-s)^{1-\alpha}>0$ for $s\in [0,t^*)$ it follows that
$$
g(s,x(s),x(\rho(s)))=0\,\,\,\,a.e. \, s\in [0,t^*).
$$
This contradicts the fact that $g(t,x,y)\neq 0$ for any $(t,x,y)\in [0,1]\times\mathbb{R}\times\mathbb{R}$. Therefore, $x(t)\neq 0$ for $t\in (0,1)$.

As an application of Bolzano's theorem, we have the following corollary.

\textbf{Corollary 1.} Under assumptions of Proposition 2, the solution $x(t)$ of Problem (1) obtained in Theorem 3 satisfies that $x(t)>0$ for $t\in (0,1)$ or $x(t)<0$ for $t\in (0,1)$.

\bigskip

On the other hand, if we perturb the data function in Problem (1) of the following way
$$
\left\{
\begin{aligned}
&D_{0+}^\alpha \left[\frac{x(t)}{f(t,x(t),x(\varphi(t)))}\right]=g(t,x(t), x(\rho(t)))+\eta(t),\,0<t<1\\
&x(0)=0,
\end{aligned}
\right.
$$
where $\eta\in C[0,1]$, $\alpha\in (0,1),\varphi,\rho\in C[0,1]$, $f\in C\left([0,1]\times\mathbb{R}\times\mathbb{R},\mathbb{R}\setminus \{0\}\right)$, and
$g\in C\left([0,1]\times\mathbb{R}\times\mathbb{R},\mathbb{R}\right)$. In this case, assumptions (H1), (H2) and (H3) of Theorem 3 if $f(t,x,y)$ and $g(t,x,y)$ satisfy (H3) of Theorem 3 and, only, we would have to check assumption (H4). This fact makes that Theorem 2 is applicable to a great number of cases.

In the sequel, we compare our results with ones appearing in the literature. In [13], the authors studied the fractional hybrid differential equation
\begin{equation}\label{eq3}
\left\{
\begin{aligned}
&D_{0+}^\alpha \left[\frac{x(t)}{f(t,x(t))}\right]=g(t,x(t)),\,\,a.e.\, t\in [0,1]\\
&x(0)=0
\end{aligned}
\right.
 \end{equation}
under the following conditions:

\textbf{(i)} $\alpha\in(0,1),\,f\in C\left([0,1]\times\mathbb{R}\times\mathbb{R},\mathbb{R}\setminus \{0\}\right)$ and $g\in C\left([0,1]\times\mathbb{R}\times\mathbb{R},\mathbb{R}\right)$.

\textbf{(ii)} The function $x\rightarrow\frac{x}{f(t,x)}$ is increasing in $\mathbb{R}$ almost everywhere for $t\in [0,1]$.

\textbf{(iii)} There exists a constant $L>0$ such that
$$
|f(t,x)-f(t,y)|\leq L|x-y|,
$$
for any $t\in [0,1]$ and $x,y,\in \mathbb{R}$.

\textbf{(iv)} There exists a function $h\in L^1\left([0,1],\mathbb{R}_+\right)$ such that
$|g(t,x)|\leq h(t)$ a.e. $t\in [0,1]$ for all $x\in\mathbb{R}$.

\textbf{(v)}
$$
\frac{L ||h||_{L^1}}{\Gamma(\alpha+1)}<1.
$$

In the above-mentioned paper, the authors proved the following theorem.

\textbf{Theorem 4.} Under assumptions (i)-(v), Problem (3) has a solution in $C[0,1]$.

\bigskip

It is clear that Problem (3) is a particular case of Problem (1) where the functions $f(t,x,y)$ and $g(t,x,y)$ are independent of $y$.

Notice that if in Example 1 we consider as $f(t,x,y)=\frac{1}{\alpha}\sqrt[4]{1+|x|}$ and $g(t,x,y)=\frac{1}{\beta}\sqrt[3]{1+|x|}$, then this variant of Example 1 can be studied by Theorem 3 while it cannot be treated by Theorem 4 since the function $g(t,x)=\frac{1}{\beta}\sqrt[3]{1+|x|}$ does not satisfy assumption (iv) of Theorem 4.

\section{Acknowledgement}
This joint work was done during the visit of the second author (Karimov E.T.) at ULPGC, supported by ERASMUSMUNDUS TIMUR project (Academic Staff mobility). The third author was partially supported by the project MTM2013-44357-P.


\textit{Bibliography.}

\end{document}